\documentclass[12pt]{amsart}
\usepackage{amssymb}
\usepackage{eucal}
\newcommand{\rr}{\mathbb{R}}
\newcommand{\p}{\partial}

\newtheorem{corollary}{Corollary}[section]
\newtheorem{lemma}{Lemma}[section]

\newtheorem{theorem}{Theorem}[section]

\author{L. Escauriaza, C. E. Kenig, G. Ponce, and L. Vega}
\title[Schr\"odinger Equations]
{On unique continuation of solutions of 
Schr\"odinger  equations}

\begin{document}
\maketitle

\numberwithin{equation}{section}

\section{Introduction}

 In this paper we study unique continuation  properties of solutions of Schr\"odinger
equations. In this class we shall include  linear ones of the form
\begin{equation}
\label{1.1}
i\p_t u + \Delta u =V u,\;\;\;\;(x,t)\in\Bbb R^n\times \Bbb R,
\end{equation}
and non-linear ones of the type
\begin{equation}
\label{1.2}
i\p_t u + \Delta u + F(u,\overline u)=0,\;\;\;\;(x,t)\in\Bbb R^n\times \Bbb R.
\end{equation}
 
 Our goal is to obtain sufficient conditions on the behavior of the solution $u$ at two different times 
 $t_0=0$ and $t_1=1$ which guarantee that $u\equiv 0$ is the unique solution
 of \eqref{1.1}. In the case of the nonlinear equation \eqref{1.2} we are interested in deducing 
 uniqueness of the solution from  information on the 
 difference of two possible solutions at  two different times.
 
   For dispersive models, and in particular  
 for Schr\"odinger equations, these kind of uniqueness results have been obtained under the assumptions that
 the solutions coincide in a large sub-domain of $\Bbb R^n$ at  two different times.
 
  For the $1$-D cubic Schr\"odinger equation , i.e. $F=\pm|u|^2 u$, $n=1$ in \eqref{1.2},
  B.-Y. Zhang \cite{Zh2} showed that if $u=0$ in $(-\infty,a)\times \{t_0, t_1\}$ 
  (or in $(a,\infty)\times \{t_0, t_1\}$)
  for some $a\in \Bbb R$, then $u\equiv 0$. His proof  is based on inverse scattering theory (IST).
  
  In \cite{KePoVe2} it was shown, under general assumptions on $F$ in \eqref{1.2}, that  if 
   $u_1, u_2\in C([0,1]:H^s(\Bbb R^n)),$
   $ s>n/2,$ $\;s\geq 2,$ such that
   \begin{equation}
\label{1.3}
   u_1(x,t)=u_2(x,t),\;\;\;\;(x,t)\in\Gamma_{x_0}^c\times \{t_0, t_1\},
 \end{equation}
   where $\Gamma_{x_0}^c$ denotes the complement of a cone $\Gamma_{x_0}$ with vertex  in $x_0\in \Bbb R^n$ with opening $< 180^0$,
   then $u_1\equiv u_2$. In \cite{KePoVe2} one of the key steps in the proof was a uniform exponential
   decay estimate in the time interval $[0,1]$ obtained under the assumption that the corresponding solution has the same decay at times $t_0=0, \,t_1=1$ (see Lemma 2.1 below). The proof of this estimate follows by combining energy estimates
  for the Fourier transform of the solution and its projection onto the positive and negative frequencies 
  together with some classical estimates for pseudo-differential and singular integral operators. 
  Roughly, this uniform estimate allows to extend (1.3) to the whole time interval $[0,1]$.
 In  this setting  one can then use V. Isakov's approach  in \cite{Is} to obtain the desired result.
  
  In \cite{IoKe1}-\cite{IoKe2}, the results in \cite{KePoVe2} were extended to the case of semi-spaces, i.e. cones with opening $=180^0$, to more general classes of potentials, and to less regular solutions.

  Unique continuation of the kind described above has been
  also established in other dispersive equations. In particular, L. Robbiano \cite{Ro} proved
 the following uniqueness result for the Korteweg-de Vries (KdV) equation.
 If two solutions $u_1,\,u_2$ of this equation, in an appropriate class, agree in a semi-line $(a,\infty)$ at time $t_0=0$ and 
for any $(x,t)\in (a,\infty)\times[0,1]$ their difference $u_1-u_2$ and its space derivatives up to order  $2$  are point-wise bounded by
$c\,e^{-x^{\alpha}}$, for some $\alpha>9/4$, then $u_1\equiv u_2$. 

Similar uniqueness result for the KdV but under the assumption that one of the solutions is $u(x,t)\equiv 0$ and the other vanishes in a semi-line
at the  times $0,\,1$ were previously proved in \cite{Zh1} by using the IST. 
This result was extended in \cite{KePoVe1} to any pair of solutions to  the generalized KdV equation, which includes
non-integrable models. More recently, based on the IST,  S. Tarama  \cite{Ta} proved that if the initial data has an appropriated exponential decay  for $x>0$, then the corresponding solution
of the KdV becomes analytic respect  the $x$ variable for all $t>0$.  We notice that even in the KdV case neither of the results  in \cite{Ro} and \cite{Ta} described above implies the other one, since in \cite{Ro}
the decay assumption is needed in a whole time interval $[0,1]$, and the result in \cite{Ta} does not apply to the difference of two arbitrary solutions of the KdV. 

Our motivation came in part from the following result due to G. H. Hardy (see \cite{StSh}) concerning 
the decay of a function and its Fourier transform. For $f:\Bbb R\to\Bbb C$ if 
$f(x)=O(e^{-\pi Ax^2})$ and its Fourier transform $\hat f(\xi)=O(e^{-\pi B\xi^2})$ with
$A, B>0$ and $AB>1$, then $f\equiv 0$. Also, if $A=B=1$, then $f(x)=ce^{-\pi x^2}$.

This kind of uncertainty principle can be re-phased in terms of the solution, 
$v(x,t)=e^{it\p_x^2/4\pi}v_0(x)$,
of the free Schr\"odinger equation
\begin{equation}
\label{1.4}
i\p_t v + \frac{1}{ 4 \pi}\p_x^2 v=0,
\end{equation}
since 
$e^{\pi i |x|^2}v(x,1)=e^{\pi i |x|^2}e^{i\partial^2_x/4\pi  }v_0(x)$ is the Fourier transform of 
$e^{i\pi\xi^2}\hat v_0(\xi)$.

Our main result concerning the equation \eqref{1.1} is the following.

 \begin{theorem}
 \label{Theorem 1.1}
 
 Let $u\in C([0,1]:H^2(\Bbb R^n))$ be a strong solution of the equation (1.1) in the domain
 $(x,t)\in\Bbb R^n\times [0,1]$ with $V:\Bbb R^n\times [0,1]\to\Bbb C$, $V\in L^{\infty}(\Bbb R^n_x
 \times[0,1])$, and
 $\nabla_x V\in L^1_t([0,1]:L^{\infty}(\Bbb R^n))$.  If there exist $\alpha>2$ and $a>0$ such that
 \begin{equation}
\label{1.5}
 u_0=u(\cdot,0),\;\;\;u_1=u(\cdot,1)\in H^1(e^{a|x|^{\alpha}}dx),
 \end{equation}
 and
 \begin{equation}
\label{1.6}
\lim_{r\uparrow \infty} \,\| V\|_{L^1_tL^{\infty}_x\{|x|>r\}}
= \lim_{r\uparrow \infty}\,\int_0^1\,\sup_{|x|>r}\,|V(x,t)|\,dt=0,
 \end{equation}
 then $u\equiv 0$. 
 
 \end{theorem}
 
  For the nonlinear equation \eqref{1.2} we shall prove,
  
 \begin{theorem}
 \label{Theorem 1.2}
 
 Let $u_1,\,u_2\in C([0,1]:H^k(\Bbb R^n)),\,k\in \Bbb Z^+,\,k>n/2+1$ be strong solutions of the equation (1.2) in the domain
 $(x,t)\in\Bbb R^n\times [0,1]$, with $F:\Bbb C^2\to \Bbb C$, $F\in C^k$  and $F(0)=\partial_uF(0)=\partial_{\bar u}F(0)=0$.
 
 If there exist $\alpha>2$ and $a>0$  such that
 \begin{equation}
\label{1.7}
  w_0=u_1(\cdot,0)-u_2(\cdot,0),\;\;w_1=u_1(\cdot,1)-u_2(\cdot,1)\in H^1(e^{a|x|^{\alpha}}dx),
 \end{equation}
 then $u_1\equiv u_2$. 
 
 \end{theorem}
 
 We shall say that  $f\in H^1(e^{a|x|^{\alpha}}dx)$ if $f,\,\partial_{x_j} f\in L^2(e^{a|x|^{\alpha}}dx)$ for
 $j=1,..,n$, i. e. 
 \begin{equation}
 \label{1.8}
 \int_{\Bbb R^n}\,|f(x)|^2\,e^{a|x|^{\alpha}}dx + \,\sum_{j=1}^n \,\int_{\Bbb R^n}\,|\partial_{x_j} f(x)|^2\,e^{a|x|^{\alpha}}dx< \infty.
 \end{equation}

 \underbar{Remarks}
 
 a) It will be clear from our proof below that Theorems 1.1 and 1.2 still hold assuming  \eqref{1.5} and \eqref{1.7} respectively with $\alpha=2$ and $a\geq c_0$, where 
 $$
 c_0=c_0(\|u\|_{L^{\infty}H^2};n;\|V\|_{L_{t,x}^{\infty}};\|\nabla_xV\|_{L^1_tL^{\infty}_x})>0,\;\;\text{ in Theorem 1.1},
 $$
 and
 $$
 c_0=c_0(\|u\|_{L^{\infty}H^2};n;\|F\|_{C^k})>0,\;\;\;\text{ in Theorem 1.2}.
 $$

 b) Theorems \ref{Theorem 1.1} and \ref{Theorem 1.2}  can be seen as  natural  extensions to the Schr\"odinger equation of the results  we recently obtained in \cite{EsKePoVe}
 for the heat equation. In \cite{EsKePoVe},  the decay assumption was only
 assumed at time $t=1$.

 c) The method of proof of  Theorems \ref{Theorem 1.1} and \ref{Theorem 1.2} follows a similar argument, which is based on two main steps. The first one is based on the exponential
 decay estimates obtained in \cite{KePoVe2}. These estimates are  expressed in terms of the $L^2$-norm with respect to the  measure
 $e^{\beta|x|}dx $ and involve bounds independent of $\beta$. Here we shall use them to deduce similar ones 
 but with higher order power in the exponent, i.e. with super-linear growth. The second main step  is to establish 
 asymptotic lower bounds for the $L^2$-norm of the solution and its space gradient  on the annulus domain 
 $(x,t)\in \{R-1<|x|<R\}\times[0,1]$. The use of a lower bound was motivated by the recent work of J. Bourgain and C. E. Kenig \cite{BoKe} on a class of stationary Schr\"odinger operators ( i.e.  $-\Delta+V(x)$).  There  a key lower bound on the decay of the average of the solution over unit balls was deduced in terms of the size of their centers. In this second part we also follow some of the arguments found in \cite{Is}.
 
 d) In a forthcoming work we shall extend the results obtained here to the generalized KdV equation.

 e) For the existence of solutions of the IVP associated to the equations \eqref{1.1}  and \eqref{1.2}
 with data at $t=0$ as in Theorem 1.1 and 1.2 we refer to \cite{Ca} and references therein.

 The rest of this paper is organized as follows. In section 2, we deduce the energy estimate
 with super-linear exponential growth in the interval $[0,1]$ under the assumption that a similar one holds
 for the solution at times $t_0=0$ and $t_1=1$. In section 3, we shall obtain a lower bound
 for the $L^2$-norm of the solution and its gradient in  the annular domain mentioned above.
 Finally, in section 4 we combine the results in the previous  sections to prove Theorems \ref{Theorem 1.1} and \ref{Theorem 1.2} .

\section{weighted energy estimates}

 In \cite{KePoVe2}  the following exponential decay estimate was established.

\begin{lemma}
\label{Lemma 2.1}

 There exists $\epsilon>0$ such that if $ V:\Bbb R^n\times [0,1]\to\Bbb C$ satisfies that
\begin{equation}
\label{2.1}
\|V\|_{L^1_tL^{\infty}_x}\leq \epsilon,
\end{equation}
and $u\in C([0,1]:L^2_x(\Bbb R^n))$ is a (strong) solution of the IVP
\begin{equation}
\label{2.2}
\begin{cases}
\begin{aligned}
&i\p_t u + \Delta u = V u + H,\;\;\;\;\;(x,t)\in \Bbb R^n\times [0,1],\\
&u(x,0)=u_0(x).
\end{aligned}
\end{cases}
\end{equation}
with $H\in L^1_t([0,1]:L^2_x(e^{2\beta x_1}dx))$, and for some $\beta\in\Bbb R$
\begin{equation}
\label{2.3}
u_0,\,u_1\equiv u(\cdot,1)\in L^2(e^{2\beta x_1}dx),
\end{equation}
then
\begin{equation}
\label{2.4}
\begin{aligned}
&\sup_{0\leq t\leq 1}\|u(\cdot,t)\|_{L^2(e^{2\beta x_1}dx)}\\
\\
&\leq c (\|u_0\|_{L^2(e^{2\beta x_1}dx)} + \|u_1\|_{L^2(e^{2\beta x_1}dx)}
+\|H\|_{L^1_tL^2_x(e^{2\beta x_1}dx)}),
\end{aligned}
\end{equation}
with $c$ independent of $\beta$.
\end{lemma}

 We shall use Lemma 2.1 to deduce further weighted inequalities.

\begin{corollary}
\label{Corollary 2.2}

  If in addition to the hypothesis in Lemma 2.1 one has that for some $\beta>0$
  \begin{equation}
\label{2.5}
  u_0,\,u_1\in L^2(e^{2\beta |x|}dx),
  \end{equation}
  and $H\in L^1_t([0,1]:L^2_x(e^{2\beta |x|}dx))$, then
 \begin{equation}
\label{2.6}
  \begin{aligned}
&\sup_{0\leq t\leq 1}\|u(\cdot,t)\|_{L^2(e^{2\beta |x|/\sqrt{n}}dx)}\\
\\
&\leq c (\|u_0\|_{L^2(e^{2\beta |x|}dx)} + \|u_1\|_{L^2(e^{2\beta |x|}dx)}
+\|H\|_{L^1_tL^2_x(e^{2\beta |x|}dx)}),
\end{aligned}
\end{equation}
 with $c$ independent of $\beta>0$.

\end{corollary}

\begin{proof}
\label{Proof of Corollary 2.2}  

 From (2.4) it follows that for any $\beta>0$ 
 \begin{equation}
\label{2.7}
\begin{aligned}
&\sup_{0\leq t\leq 1}\|u(\cdot,t)\|_{L^2(e^{\pm 2\beta x_j}dx)}\\
\\
&\leq c (\|u_0\|_{L^2(e^{2\beta|x|}dx)} + \|u_1\|_{L^2(e^{2\beta|x|}dx)}
+\|H\|_{L^1_tL^2_x(e^{2\beta|x|}dx)})\equiv \varPhi,
\end{aligned}
\end{equation}
for any $j=1,..,n$. Hence, for any $t\in[0,1]$ one has that
\begin{equation}
\label{2.8}
\begin{aligned}
&\|u(\cdot,t)\|_{L^2(e^{2\beta |x|/\sqrt{n}}dx)}
\leq \sum_{j=1}^n(\int_{\Bbb R^n}|u(x,t)|^2\,e^{2\beta|x_j|}\,dx)^{1/2}\leq 2n\varPhi.
\end{aligned}
\end{equation}
Taking the supremum on $t\in[0,1]$ we get the desired inequality (2.6).
\end{proof} 

\vskip.1in

\begin{corollary}
\label {Corollary 2.3}

 If in addition to the hypothesis in Lemma 2.1 one assumes that for some $a>0$ and $\alpha>1$
  \begin{equation}
\label{2.9}
  u_0,\,u_1\in L^2(e^{a |x|^{\alpha}}dx),
 \end{equation}
  and $H\in L^1_t([0,1]:L^2_x(e^{a |x|^{\alpha}}dx))$, with $\,u\in C([0,1]:H^1(\Bbb R^n))$, then 
  there exists $c_{\alpha}>0$ such that 
 \begin{equation}
\label{2.10}
 \begin{aligned}
&\sup_{0\leq t\leq 1}\int_{|x|\geq c_{\alpha}}\,|u(x,t)|^2\,e^{a|x|^{\alpha}/(10 \sqrt{n})^{\alpha}}dx\\
\\
&\leq c (\|u_0\|^2_{L^2(e^{a |x|^{\alpha}}dx)} + \|u_1\|^2_{L^2(e^{a|x|^{\alpha}}dx)})\\
&+c\,\int_0^1\,\int_{\Bbb R^n}|H(x,t)|^2e^{a|x|^{\alpha}}dx\,dt+c \sum_{l=0}^1\int_0^1
\int_{\Bbb R^n} |\partial^l_r u(x,t)|^2\,dx\,dt.
\end{aligned}
\end{equation}

The constant $c$ in (2.10) can be taken uniform in a set of $a$'s bounded from below away from zero.

\end{corollary}

\begin{proof}
\label{Proof of Corollary 2.3}

 We multiply the equation  in (2.2) by $\eta_R(x)=\eta(x/R)$ with $\eta\in C^{\infty}$
 non-decreasing radial function such that $\eta(x)=0$ if $|x|<1$ and $\eta(x)=1$ if $|x|>2$.
 
 Using the notation $u_R(x,t)=u(x,t)\,\eta_R(x)$ we get the new equation
 \begin{equation}
\label{2.11}
 i\partial_t u_R+\Delta u_R=Vu_R+\tilde H_R,
 \end{equation}

 with
 \begin{equation}
\label{2.12}
 \tilde H_R=H\eta_R-2\partial_r\eta_R\partial_ru-(\partial_r^2\eta_R+\frac{n-1}{r}\partial_r
 \eta_R)u
 \end{equation}
 Applying  Corollary 2.1 (estimate \eqref{2.6}) it follows that
\begin{equation}
\label{2.13}
 \begin{aligned}
 &\int_{|x|>2R}\,|u(x,t)|^2\,e^{2\beta|x|/\sqrt{n}}dx
 \leq c\,\sum_{j=0}^1\,\int_{|x|>R} |u_j(x)|^2\,e^{2\beta|x|}dx\\
 &+ c\,\int_0^1\,\int_{|x|>R}|H(x,t)|^2e^{2\beta|x|}dx dt\\
 &+c\int_0^1 \,\int_{R<|x|<2R}
 \,(\frac{|\partial_ru|^2}{R^2}
 +\frac{|u|^2}{R^4})\,e^{2\beta|x|}dx dt.
 \end{aligned}
 \end{equation}
 Since
 \begin{equation}
\label{2.14}
 \begin{aligned}
 &\int_0^1 \,\int_{R<|x|<2R}
 \,(\frac{|\partial_ru|^2}{R^2}
 +\frac{|u|^2}{R^4})\,e^{2\beta|x|}dx dt\\
 &\leq \,e^{4\beta R}
 \,\int_0^1 \,\int_{R<|x|<2R}
 \,(\frac{|\partial_ru|^2}{R^2}
 +\frac{|u|^2}{R^4})dx dt,
 \end{aligned}
 \end{equation}
multiplying the expression in \eqref{2.13}
by $\,e^{-4\beta R}$  one has  that
\begin{equation}
\label{2.15}
 \begin{aligned}
 &A_1\equiv e^{-4\beta R}\,\int_{|x|>2R}\,|u(x,t)|^2\,e^{2\beta|x|/\sqrt{n}}dx\\
 &
 \leq c\,e^{-4\beta R}\,\sum_{j=0}^1\,\int_{|x|>R} |u_j(x)|^2\,e^{2\beta|x|}dx\\
 &+ c\,e^{-4\beta R}\,\int_0^1\,\int_{|x|>R}|H(x,t)|^2e^{2\beta|x|}dx dt\\
 &+c \,\int_0^1 \,\int_{R<|x|<2R}
 \,(\frac{|\partial_ru|^2}{R^2}
 +\frac{|u|^2}{R^4})dx dt\equiv B_1+B_2+B_3.
 \end{aligned}
 \end{equation}

Next, we fix $2\beta=bR^{\alpha-1}$, with $b=b(\alpha)>0$ to be determined,
integrate the inequality \eqref{2.15} in $R$ in the interval $[1,\infty)$, and consider 
the resulting terms separately. First, for the term coming from $B_1$
using Fubini's theorem we write
\begin{equation}
\label{2.16}
\begin{aligned}
&\int_1^\infty e^{-2b R^{\alpha}}\,\sum_{j=0}^1\,\int_{|x|>R} |u_j(x)|^2\,e^{2\beta|x|}dx \,dR\\
\\
&=\sum_{j=0}^1 \,\int_{|x|>1}(\int_1^r \,e^{-2bR^{\alpha}+bR^{\alpha-1}r}dR)\,
 |u_j(x)|^2\,dx,
\end{aligned}
\end{equation}
where $r=|x|$. To deduce an upper bound for this expression we see that 
$\varphi(R)=bR^{\alpha-1}(r-2R)$ has its 
maximum at $R_M=(\alpha-1)r/2\alpha<r/2$, hence
\begin{equation}
\label{2.17}
\begin{aligned}
&\int_1^r \,e^{-2bR^{\alpha}+bR^{\alpha-1}r}dR
\leq re^{\varphi(R_M)}\\
&=r\,e^{b(\alpha-1)^{\alpha-1}r^{\alpha}/(2^{\alpha-1}\alpha^{\alpha-1})}
=r\,e^{b_{\alpha}r^{\alpha}}=|x|\,e^{b_{\alpha}|x|^{\alpha}},
\end{aligned}
\end{equation}
i.e. $b_{\alpha}=b(\alpha-1)^{\alpha-1}/(2^{\alpha-1}\alpha^{\alpha})$.
This estimated inserted in \eqref{2.16} yields the bound
\begin{equation}
\label{2.18}
\sum_{j=0}^1\,\int_{\Bbb R^n} |u_j(x)|^2 e^{b_{\alpha}|x|^{\alpha}}\,|x|\,dx.
\end{equation}
A similar argument provides the following upper bound for the term coming from $B_2$ 
in \eqref{2.15}
\begin{equation}
\label{2.19}
\int_0^1\,\int_{\Bbb R^n} |H(x,t)|^2e^{b_{\alpha}|x|^{\alpha}}\,|x|\,dx dt.
\end{equation}

Next, we shall deduce a lower bound for the term arising from $A_1$. Using again Fubini's theorem
 this can be written as
\begin{equation}
\label{2.20}
 \begin{aligned}
&\int_1^{\infty} e^{-2b R^{\alpha}}\,\int_{|x|>2R}\,|u(x,t)|^2\,e^{2\beta|x|/\sqrt{n}}dx dR\\
&=\int_2^{\infty}	 \int_{\Bbb S^{n-1}}(\int_1^{r/2} \,e^{-2bR^{\alpha}+bR^{\alpha-1}r/\sqrt{n}}dR)
|u(x,t)|^2\,r^{n-1}dS\,dr.
\end{aligned}
\end{equation}
Since $\,\eta(R)= -2bR^{\alpha}+bR^{\alpha-1}r/\sqrt{n}\,$ has its maximum at 
$\tilde R_M=
(\alpha-1)r/2\alpha\sqrt{n}<r/2\sqrt{n}< r/2$, we take $R_0=(\alpha-1)r/10\alpha\sqrt{n}$,
$r>c_{\alpha}$ and $c_{\alpha}>10\alpha \sqrt{n}/(\alpha-1)>2$ to bound
from below
the integral inside \eqref{2.20}  as 
\begin{equation}
\label{2.21}
\begin{aligned}
&\int_1^{r/2} \,e^{-2bR^{\alpha}+bR^{\alpha-1}r/\sqrt{n}}dR
>\int_{R_0}^{\tilde R_M}\,e^{bR^{\alpha-1}(r/\sqrt{n}-2R)}dR\\
&>e^{bR_0^{\alpha-1}(r/\sqrt{n}-2\tilde R_M)}(\tilde R_M-R_0)\\
&\geq\frac{2}{5}\frac{\alpha-1}{\alpha}\frac{r}{\sqrt{n}}\,e^{b(\alpha-1)^{\alpha-1}r^{\alpha}/(10^{\alpha-1}\alpha^{\alpha}\sqrt{n}^{\alpha})}\\
&=
\frac{2}{5}\frac{\alpha-1}{\alpha}\frac{r}{\sqrt{n}}\,e^{b_{\alpha}r^{\alpha}/(5^{\alpha-1}  \alpha
\sqrt{n}^\alpha)}.
\end{aligned}
\end{equation}
The term obtained from $B_3$ in \eqref{2.15} can be handled as 
\begin{equation}
\label{2.22}
\begin{aligned}
&\int_1^{\infty}
\int_0^1 \,\int_{R<|x|<2R}
 \,\left(\frac{|\partial_ru|^2}{R^2}+\frac{|u|^2}{R^4}\right)dx dt dR\\
 &c \leq \int_0^1\,\int_{|x|\geq 1}
\,(|\partial_ru|^2
 +|u|^2) (\int_{|x|/2}^{|x|} \frac{dR}{R^2+R^4})dx dt \\
& \leq c \int_0^1\,\int_{\Bbb R^n} (|u(x,t)|^2+|\partial_ru(x,t)|^2)dx dt.
 \end{aligned}
 \end{equation}
 
  Thus, collecting this information and those
 in the estimates \eqref{2.15}-\eqref{2.22}, fixing $b$ such that $\,b_{\alpha}+\epsilon=b(\alpha-1)^{\alpha-1}/(2^{\alpha-1}\alpha^{\alpha})+\epsilon=a$,
 with $\,\epsilon>0$ small enough we obtain the desired estimate \eqref{2.10}.

\end{proof}

  A similar argument provides the following result.

\begin{corollary}
\label{Corollary 2.4}

 If in addition to the hypothesis in Lemma 2.1 one assumes that for some $\,a\in \Bbb R$ and $\alpha>1$
  \begin{equation}
\label{2.23}
  u_0,\,u_1\in L^2(e^{a x_1 |x_1|^{\alpha-1}}dx),
 \end{equation}
  and $H\in L^1_t([0,1]:L^2_x(e^{a x_1 |x_1|^{\alpha-1}}dx))$, with  $\,u\in C([0,1]:H^1(\Bbb R^n))$,  then 
  there exists $c_{\alpha}$ such that
 \begin{equation}
\label{2.24}
 \begin{aligned}
&\sup_{0\leq t\leq 1}\int_{|x_1|\geq c_{\alpha}}\,|u(x,t)|^2\,e^{a x_1 |x_1|^{\alpha-1}/(10)^{\alpha}}dx\\
\\
&\leq c (\|u_0\|^2_{L^2(e^{a x_1 |x_1|^{\alpha-1}}dx)} + \|u_1\|^2_{L^2(e^{a x_1 |x_1|^{\alpha-1}}dx)})\\
&+c\,\int_0^1\,\int_{\Bbb R^n}|H(x,t)|^2e^{a x_1 |x_1|^{\alpha-1}}dx\,dt\\
&+c \sum_{l=0}^1\int_0^1
\int_{\Bbb R^n} |\partial^l_{x_1} u(x,t)|^2\,dx\,dt.
\end{aligned}
\end{equation}
   
  \end{corollary}
  
 In Corollaries 2.2 and 2.3 it suffices to assume that $u,\, \partial_r u$ and $u,\,\partial_{x_1} u\,$ belong to $\,C([0,1]:L^2(\Bbb R^n))$, respectively.

Also, the results in this section extend to equations of the form
\begin{equation}
\label{2.25}
i\partial_t u +\Delta u = V_1 u + V_2 \bar u + H,
\end{equation}
with the potentials $V_j(x,t),\;j=1,2$ satisfying the assumption (2.1).

\section{lower bounds estimates}

 This section is concerned with lower bounds  for the $L^2$-norm of the
 solution of the equations \eqref{1.1} and  \eqref{1.2} and its  gradient in the domain 
  $\{R-1<|x|<R\}\times[0,1]$ . 
 
 \begin{lemma}
 \label{Lemma 3.1}
Assume that $R>0$ and $\varphi:[0,1]\longrightarrow\rr$ is a smooth function. Then, there exists $c=c(n, \|\varphi'\|_{\infty}+\|\varphi''\|_{\infty})>0$ such that, the inequality
\begin{equation}
\label{3.1}
\frac{\alpha^{3/2}}{R^2}\, \|e^{\alpha |\frac xR+\varphi(t)e_1|^2}g\|_{L^2(dxdt)}\leq
c\| e^{\alpha |\frac xR+\varphi(t)e_1|^2}(i \partial_t+\Delta)g\|_{L^2(dxdt)}
\end{equation}
holds, when $\alpha \ge cR^2$ and $g\in C_0^\infty(\rr^{n+1})$ has its support contained in the set \[\{(x,t): |\tfrac xR+\varphi(t)e_1|\ge 1\}\ .\]
\end{lemma}
\begin{proof}
\label{Proof of Lemma 3.1}
We shall follow the arguments in \cite{Is} and \cite{IoKe2}. Let $f=e^{\alpha |\frac xR+\varphi(t)e_1|^2}g$.
Then,
\begin{equation}
\label{3.3}
e^{\alpha|\frac xR+\varphi(t)e_1|^2}(i \partial_t+\Delta)g
=S_\alpha f-4\alpha A_\alpha f\ ,
\end{equation}
where
\begin{align*}
&S_\alpha=i\partial_t+\Delta+\tfrac{4\alpha^2}{R^2}|\tfrac xR+\varphi e_1|^2,\\
\\
&A_\alpha=\tfrac 1R\left (\tfrac xR+\varphi e_1\right)\cdot\nabla+ \tfrac n{2R^2}+
\tfrac{i\ \varphi'}2\left(\tfrac {x_1}R+\varphi\right ).
\end{align*}
Thus,
\begin{equation}
\label{3.5}
A_\alpha^*=-A_\alpha,\;\;\;\;\;\;\;S_\alpha^*=S_\alpha,
\end{equation}
and
\begin{align*}
\|e^{\alpha |\frac xR+\varphi e_1|^2}&(i \partial_t+\Delta)g\|_2^2
=\langle S_\alpha f-4\alpha A_\alpha f, S_\alpha f-4\alpha A_\alpha f\rangle\\
&\geq -4\alpha\langle (S_\alpha A_\alpha - A_\alpha S_\alpha )f, f \rangle =  -4\alpha\langle [S_\alpha ,A_\alpha]f, f\rangle\ .
\end{align*}
A calculation shows that
\begin{equation*}
 [S_\alpha ,A_\alpha]=\tfrac 2{R^2} \Delta-\tfrac {4\alpha^2}{R^4} 
|\tfrac xR+\varphi  e_1|^2-\tfrac 12 [(\tfrac 
{x_1}R+\varphi)\varphi''+\varphi'^2]+\tfrac {2i\varphi'}R\partial_{x_1}
\end{equation*}
and
\begin{equation}
\label{3.7}
\begin{aligned}
\|e^{\alpha|\tfrac xR+\varphi e_1|^2}&(i\partial_t+\Delta)g\|_2^2\geq  \frac{16\alpha^3}{R^4} \int |\tfrac xR+\varphi e_1|^2 |f|^2dxdt+\frac {8\alpha} {R^2}\int |\nabla f|^2dxdt\\
&+2\alpha\int[(\tfrac {x_1}R+\varphi)\varphi''+\varphi'^2]|f|^2dxdt-\frac{8\alpha i}R\int \varphi'\,
\partial_{x_1}
f\bar fdxdt\ .
\end{aligned}
\end{equation}

Hence, using the hypothesis on the support on $g$ and the Cauchy-Schwarz inequality,  the absolute value of the last two terms in \eqref{3.7} can be bounded by a fraction of the  first two terms on the right hand side of \eqref{3.7}, when $\alpha\ge cR^2$ for some large $c$ depending on $\|\varphi' \|_{\infty}+\|\varphi'' \|_{\infty}$. This yields \eqref{3.1} and  Lemma \ref{Lemma 3.1}.

\end{proof}

\begin{theorem}
\label{Theorem 3.2}

Let $u\in C([0,1]:H^1(\Bbb R^n))$ be a strong solution of 
\begin{equation}
\label{3.8}
i\partial_t u +\Delta u+V u=0,\;\;\;\;t\in[0,1],\;\;x\in\Bbb R^n.
\end{equation}
If
\begin{equation}
\label{3.9}
\int_0^1\int_{\Bbb R^n}(|u|^2+|\nabla_xu|^2)(x,t)dxdt\leq A^2,
\end{equation}
\begin{equation}
\label{3.10}
\int_{1/2-1/8}^{1/2+1/8}\int_{|x|<1} |u|^2(x,t)dx dt\geq 1,
\end{equation}
and
\begin{equation}
\label{3.11}
\|V\|_{L^{\infty}(\Bbb R^n\times[0,1])}\leq L,
\end{equation}
then there exists $R_0=R_0(n,A,L)>0$ and a constant $c=c(n)$ such that  for $R\geq R_0$ 
it follows that
\begin{equation}
\label{3.12}
\delta(R)\equiv \left(\int_0^1\int_{R-1<|x|<R}(|u|^2+|\nabla_xu|^2)(x,t)dxdt\right)^{1/2}
\geq ce^{-cR^{2}}.
\end{equation}

\end{theorem}

\begin{proof}
\label{Proof of Theorem 3.2}
Let $\theta_R$, $\theta \in C^{\infty}_0(\Bbb R^n)$ and $\varphi\in C^\infty([0,1])$ be functions verifying $\theta_R(x)=1$ if $|x|\leq R-1$, $\theta_R(x)=0$ if $|x|>R$, $\theta (x)=0$  when $|x|\leq 1$, $\theta (x)=1$ if $|x|\ge 2$, $0\le \varphi \le 3$, $\varphi=3$ in the interval $[\frac 12-\frac 18,\frac 12+\frac 18]$ and $\varphi=0$ in $[0,\frac 14]\cup
 [\frac 34,1]$. 
 
 With this choice of $\varphi$, we will apply  the Lemma \ref{Lemma 3.1} to the function \[g(x,t)=\theta_R(x)\theta(\tfrac xR+\varphi(t)e_1)u(x,t)\ .\]

Observe that $g$ has compact support on $\rr^n\times (0,1)$ and  satisfies the hypothesis in Lemma \ref{Lemma 3.1}, $g=u$ in $B_{R-1}\times [\tfrac 12-\tfrac 18,\tfrac 12+\tfrac 18]$, where $|\tfrac xR+\varphi(t) e_1|\ge 3-1=2$,
\begin{equation}
\label{3.13}
\begin{aligned}
&(i\partial_t +\Delta+V)g=\theta(\tfrac xR+\varphi e_1)(2\nabla\theta_R(x)\cdot\nabla u+u\Delta 
\theta_R(x) )\\
&+\theta_R(x)\left[2R^{-1}\nabla\theta(\tfrac xR+\varphi e_1)\cdot\nabla u+R^{-2}u\Delta\theta(\tfrac xR+\varphi e_1) +i\varphi'\partial_{x_1}\theta(\tfrac xR+\varphi e_1) u \right]\ ,
\end{aligned}
\end{equation}
and that the first and second terms on the right hand side of \eqref{3.13} are supported respectively in $B_R\setminus B_{R-1}\times[0,1]$, where $|\tfrac xR+\varphi e_1|\le 4$, and in  $\{(x,t): 1\le |\tfrac xR+\varphi e_1|\le 2\}$. Thus, 
\begin{equation}
\label{3.15}
\|e^{\alpha |\frac xR+\varphi(t) e_1|^2}g\|_{L^2(dxdt)}\ge e^{4\alpha}\|u\|_{L^2(B_1\times [\frac 12-\frac 18,\frac 12+\frac 18])}\ge e^{4\alpha}\ ,
\end{equation}
and combining \eqref{3.1} and \eqref{3.11} with $\alpha\ge cR^2$ 
\begin{equation}
\label{3.16}
\frac{\alpha^{3/2}}{cR^2}\,\|e^{\alpha |\frac xR+\varphi(t)e_1|^2}g\|_{L^2(dxdt)}\le L\,\|e^{\alpha |\frac xR+\varphi(t)e_1|^2}g\|_{L^2(dxdt)}+e^{16\alpha}\delta(R)+e^{4\alpha}A\ .
\end{equation}

If we choose $\alpha=cR^2$,  it is posible to hide the first term on the right hand side of \eqref{3.16} in the left hand side of the same inequality, when $R\ge R_0(L)$. This and \eqref{3.15} imply that
\begin{equation*}
R\le c\left(e^{8cR^2}\delta(R)+A\right)\ ,
\end{equation*} 
when $R\ge R_0(L)$, which implies the Theorem \ref{Theorem 3.2} when $R\ge R_0(L,A)$.
 \end{proof}
 
\section{proof of theorems 1.1 and 1.2}

 In this section we shall combine the results in Sections 2 and 3 to prove Theorems 1.1. and 1.2.
 
 \it{Proof of Theorem 1.1}\rm
 
  If $u \not \equiv 0$, we can assume, after a possible translation, dilation,  and multiplication by a constant, that $u=u(x,t)$
  satisfies the  hypothesis of Theorem 3.1. Therefore, there exist $R_0=R_0(n,A,L)>0$ and $c=c(n)$ such  that for $R\geq R_0$ 
 \begin{equation}
 \label{4.1}
  \delta(R) =(\int_0^1\int_{R-1<|x|<R} (|u|^2+|\nabla_x u|^2)(x,t)dxdt)^{1/2}\geq c e^{-c R^{2}}.
  \end{equation}  
    
   Next, we take $\phi_R\in C^{\infty}(\Bbb R^n)$ radial with  $\phi_R(x)=0$ if $|x|<R-1$, $\phi_R(x)=1$ if
   $|x|>R$, and $\partial_r \phi_R(r)\geq 0$, to get from (1.1) the  equations
 \begin{equation}
 \label{4.2}
 \begin{aligned}
 i\partial_t (u\phi_R)+\Delta(u\phi_R)&=V_R(u\phi_R)+2\nabla_x u\cdot\nabla_x\phi_R+\Delta\phi_Ru\\
 & =V_R(u\phi_R)+H,
  \end{aligned}
  \end{equation}
   and for $j=1,..,n$
    \begin{equation}
 \label{4.3}
 \begin{aligned}
   &  i\partial_t (\partial_{x_j}(u\phi_R))+\Delta(\partial_{x_j}(u\phi_R))=V_R(\partial_{x_j}(u\phi_R))\\
   &+ \partial_{x_j}V_R(u\phi_R)
   +2\nabla_x\partial_{x_j}\phi_R\cdot\nabla_x u
   +2\nabla_x\phi_R\cdot\nabla_x\partial_{x_j}u\\
   &+\Delta\phi_R\partial_{x_j}u-\Delta\partial_{x_j}\phi_R u=
   V_R(\partial_{x_j}(u\phi_R))+\tilde H_j,
   \end{aligned}
  \end{equation}
  where $V_R(x,t)=\phi_{R-1}(x)\,V(x,t)$.
 
 In the following inequalites $c_0$ will  denote a constant independent of $R$ which may change from line to line. Now applying  Corollary 2.2 to the equation (4.2) we obtain, for $R$ large enough depending on $\alpha$, that
for $a$ large enough depending on $\alpha$, 
 \begin{equation}
 \label{4.4}
 \begin{aligned}
 &\sup_{0\leq t\leq 1}\int_{|x|>R}\,|u(x,t)|^2 e^{a|x|^{\alpha}/(10\sqrt{n})^{\alpha}}dx\\
& \leq \sup_{0\leq t\leq 1}\int_{|x|\geq c_{\alpha}}\,|(u\phi_R)(x,t)|^2 e^{a|x|^{\alpha}/(10\sqrt{n})^{\alpha}}dx\\
&\leq c\sum_{j=0}^1\|u_j\|^2_{L^2(e^{a|x|^{\alpha}}dx)}
+c\int_0^1\int_{R-1<|x|<R}(|u|^2+|\nabla_x u|^2)e^{a|x|^{\alpha}}dxdt\\
&+ c\sum_{l=0}^1\int_0^1\int_{\Bbb R^n}|\partial_r^l(u\phi_R)|^2dxdt
\leq c_0+c_0e^{aR^{\alpha}}.
 \end{aligned}
  \end{equation}
    
  Corollary 2.2 (see remark at the end of section 2)  with $a/(10\sqrt{n})^{\alpha}$ instead of $a$, the equation (4.3) and the previous estimate (4.4) leads to
 \begin{equation}
 \label{4.5}
 \begin{aligned}
 &\sup_{0\leq t\leq 1}\int_{|x|>R}\,|\nabla_xu(x,t)|^2 e^{a|x|^{\alpha}/(10\sqrt{n})^{2\alpha}}dx\\
& \leq \sup_{0\leq t\leq 1}\int_{|x|\geq c_{\alpha}}\,|\nabla_x(u\phi_R)(x,t)|^2 e^{a|x|^{2\alpha}/(10\sqrt{n})^{\alpha}}dx\\  
&\leq c_0+
c\int_0^1\int_{R-1<|x|<R}(|u|^2+|\nabla_x u|^2+|\nabla^2_x u|^2)e^{a|x|^{\alpha}}dxdt\\
&+c\int_0^1\int_{\Bbb R^n}|u\phi_R\partial_{x_j}V_R|^2
e^{aR^{\alpha} /(10\sqrt{n})^{\alpha}} dxdt\\
&\leq c_0+c_0e^{aR^{\alpha} }.
 \end{aligned}
  \end{equation}
  
 Combining (4.4)-(4.5) one sees that
 \begin{equation}
 \label{4.6}
\sup_{0\leq t\leq 1}\int_{|x|>R}\,(|u|^2+|\nabla_xu|^2)(x,t) e^{a|x|^{\alpha}/(10\sqrt{n})^{2\alpha}}dx
\leq c_0+ c_0e^{aR^{\alpha}}.
\end{equation}
 From (4.6) and (4.1) we conclude for any $\mu>1$ with $\mu R-1>R$ that 
 \begin{equation}
 \label{4.7}
 \begin{aligned}
&c_0e^{-c(\mu R)^{2}} e^{a(\mu R-1)^{\alpha}/(10\sqrt{n})^{2\alpha}}\leq c_0
\delta(\mu R) e^{a(\mu R-1)^{\alpha}/(10\sqrt{n})^{2\alpha}}\\
&\leq c_0+ c_0e^{aR^{\alpha}}.
\end{aligned}
\end{equation}
 Finally, since $\alpha>2$ taking $\mu$ sufficiently large in (4.7) we get a contradiction.
Hence, $u\equiv 0$.

\it{Proof of Theorem 1.2}\rm

 We consider the difference of the two solutions
 \begin{equation}
 \label{4.8}
 w(x,t)=u_1(x,t)-u_2(x,t),
 \end{equation}
 which satisfies the equation
 \begin{equation}
 \label{4.9}
 i\partial_t w+\Delta w=\left(\frac{F(u_1,\bar u_1)-F(u_2,\bar u_2)}{u_1-u_2}\right) w.
 \end{equation}
Also, its $\partial_{x_j}$-derivative, $j=1,..,n$, solves
\begin{equation}
 \label{4.10}
 \begin{aligned}
&i\partial_t \partial_{x_j}w+\Delta \partial_{x_j}w
=\partial_uF(u_1,\bar u_1)\partial_{x_j}w +
\partial_{\bar u}F(u_1,\bar u_1)\partial_{x_j}\bar w\\
&+\left(\frac{\partial_u F(u_1,\bar u_1)-\partial_u F(u_2,\bar u_2)}{u_1-u_2} \right) \partial_{x_j} u_2 w
\\
&+\left(\frac{\partial_{\bar u}
F(u_1,\bar u_1)-\partial_{\bar u}F(u_2,\bar u_2)}{\bar u_1-\bar u_2}\right) \partial_{x_j} \bar u_2 \bar w\\
&=\partial_uF(u_1,\bar u_1)\partial_{x_j}w +
\partial_{\bar u}F(u_1,\bar u_1)\partial_{x_j}\bar w + H_j,
\end{aligned}
\end{equation}

 The potential
\begin{equation}
 \label{4.11}
V_1(x,t)=\frac{F(u_1,\bar u_1)-F(u_2,\bar u_2)}{u_1-u_2},
\end{equation}
in the equation (4.9) satisfies the hypothesis of Theorem 3.1, therefore
\begin{equation}
 \label{4.12}
  \delta(R) =(\int_0^1\int_{R-1<|x|<R} (|w|^2+|\nabla_x w|^2)(x,t)dxdt)^{1/2}\geq c e^{-c R^{2}}.
 \end{equation}
Next, we shall follow the argument in (4.2)-(4.5). Thus, we multiply the equations in (4.9)-(4.10)
by $\phi_R(x)$ defined before (4.2) and  observe that the potentials $V_1(x,t)$ in (4.11)  and 
\begin{equation}
 \label{4.13}
V_2(x,t)= \partial_uF(u_1,\bar u_1),
\;\;\;\;V_3(x,t)=\partial_{\bar u}F(u_1,\bar u_1),
\end{equation}
satisfy that
\begin{equation}
 \label{4.14}
\|\phi_{R-1}V_j\|_{L^1_tL_x^{\infty}}\to 0,\;\;\;\text{as}\;\;R \uparrow \infty,\;\;j=1,2,3.
\end{equation}
In particular, we can apply Corollary 2.2 (see the remark at the end of Section 2)  to the equation (4.9) and argue as in (4.4)
to get that for $R$ sufficiently large 
\begin{equation}
 \label{4.15}
 \sup_{0\leq t\leq 1}\int_{|x|>R}\,|w(x,t)|^2 e^{a|x|^{\alpha}/(10\sqrt{n})^{\alpha}}dx
\leq c_0 e^{aR^{\alpha}}
\end{equation}
Also we see that the terms 
\begin{equation}
 \label{4.16}
\frac{\partial_u F(u_1,\bar u_1)-\partial_u F(u_2,\bar u_2)}{u_1-u_2} \; \partial_{x_j} u_2\,w,
\end{equation}
and
\begin{equation}
 \label{4.8}
\frac{\partial_{\bar u}
F(u_1,\bar u_1)-\partial_{\bar u}F(u_2,\bar u_2)}{\bar u_1-\bar u_2}\; \partial_{x_j} \bar u_2\,\bar w,
\end{equation}
in the equation (4.10) belong to $L^{\infty}_t([0,1]:L^{\infty}_x(\Bbb R^n))$, and contain a factor 
already estimated in (4.15). Also, the additional terms coming from the commutator between multiplication by $\phi_R$ and the Laplacian $\Delta$ are similar to those considered in 
(4.2)-(4.3), and their  bounds are analogous to those described in (4.4)-(4.5). Hence, the reminding of  the proof of Theorem 1.2 follows the argument 
given in (4.2)-(4.7) in the proof of Theorem 1.1, therefore it will be omitted.

\vspace{3mm} \noindent{\large {\bf Acknowledgments}}
\vspace{3mm}\\  L. E. and L. V. were supported by a MEC grant and by the European Comission via the network Harmonic Analysis and Related Problems. C. E. K. and G. P. were supported by NSF grants.

%\end{document}

\vskip1cm

\noindent{\bf Luis Escauriaza}\\
Departamento de Matematicas\\
Universidad del Pais Vasco\\
Apartado 644\\
48080 Bilbao\\
Spain\\
E-mail: mtpeszul@lg.ehu.es\\

\noindent{\bf Carlos E. Kenig}\\
Department of Mathematics\\
University of Chicago\\
Chicago, Il. 60637 \\
USA\\
E-mail: cek@math.uchicago.edu\\

\noindent{\bf Gustavo Ponce}\\
Department of Mathematics\\
University of California\\
Santa Barbara, CA 93106\\
USA\\
E-mail: ponce@math.ucsb.edu\\

\noindent{\bf Luis Vega}\\
Departamento de Matematicas\\
Universidad del Pais Vasco\\
Apartado 644\\
48080 Bilbao\\
Spain\\
E-mail: mtpvegol@lg.ehu.es

\end{document}